\documentclass[12pt]{article}
\usepackage{amssymb}
\usepackage{latexsym}
\usepackage{amsmath}
\usepackage{amsthm}
\usepackage{comment}

\def\e{\varepsilon}
\def\P{\mathbb P}
\def\E{\mathbb E}
\def\R{\mathbb R}
\newtheorem{lem}{Lemma}
\newtheorem{thm}{Theorem}

\theoremstyle{definition}

\usepackage{xpatch}
\xpatchcmd{\proof}{\itshape}{\normalfont\proofnameformat}{}{}
\newcommand{\proofnameformat}{}

\begin{document}

\tolerance=1000

\renewcommand{\proofnameformat}{\bfseries}

\title{\bf On the discrepancy of random subsequences of $\{n\alpha\}$ II}
\author{Istv\'an Berkes\footnote{Alfr\'ed R\'enyi Institute of Mathematics, 1053 Budapest, Re\'altanoda u.\ 13-15, Hungary. e-mail: \texttt{berkes.istvan@renyi.hu}. Research supported by NKFIH grant K 125569.} \, and Bence Borda\footnote{Alfr\'ed R\'enyi Institute of Mathematics, 1053 Budapest, Re\'altanoda u.\ 13-15, Hungary and Graz University of Technology, 8010 Graz, Steyrergasse 30, Austria. e-mail: \texttt{borda.bence@renyi.hu}. Research supported by the Austrian Science Fund (FWF), project Y-901.
}}
\date{}
\maketitle

\abstract{Let $\alpha$ be an irrational number, let $X_1, X_2, \ldots$ be independent, identically distributed, integer-valued random variables, and put $S_k=\sum_{j=1}^k X_j$. Assuming that $X_1$ has finite variance or heavy tails $\P (|X_1|>t)\sim ct^{-\beta}$, $0<\beta<2$, in \cite{BB2} we proved that, up to logarithmic factors, the order of magnitude of the discrepancy $D_N (S_k \alpha)$ of the first $N$ terms of the sequence $\{S_k \alpha\}$ is $O(N^{-\tau})$, where $\tau= \min (1/(\beta \gamma), 1/2)$ (with $\beta=2$ in the case of finite variances) and $\gamma$ is the strong Diophantine type of $\alpha$. This shows a change of  behavior of the discrepancy at $\beta\gamma=2$. In this paper we determine the exact order of magnitude of $D_N (S_k \alpha)$ for $\beta\gamma<1$, and determine the limit distribution of $N^{-1/2} D_N (S_k \alpha)$. We also prove a functional version of these results describing the asymptotic behavior of a wide class of functionals of the sequence $\{S_k \alpha\}$. Finally, we extend our results to the discrepancy of $\{S_k\}$ for general random walks $S_k$ without arithmetic conditions on $X_1$, assuming only a mild polynomial rate on the weak  convergence of $\{S_k\}$ to the uniform distribution.}

\vspace{5mm}

\noindent\textbf{MSC2020:} 11K38, 11K60, 60G50, 60F17

\vspace{5mm}

\noindent\textbf{Keywords:} i.i.d.\ sums, random walk, discrepancy, Diophantine approximation, functional central limit theorem, functional law of the iterated logarithm

\newpage
\section{Introduction}

Let $\alpha$ be an irrational number with strong Diophantine type $\gamma$; that is,
\[ 0< \liminf_{q \to \infty} \| q \alpha \| q^{\gamma} < \infty , \]
where $\| \cdot \|$ denotes distance from the nearest integer. Let $X_1, X_2, \ldots$ be independent, identically distributed (i.i.d.) integer-valued random variables, and put $S_k= \sum_{j=1}^k X_j$. Assuming that $X_1$ has finite variance or heavy tails $\P (|X_1|>t)\sim ct^{-\beta}$, $0<\beta<2$,
in our previous paper \cite{BB2} we investigated the order of magnitude of the discrepancy $D_N(S_k \alpha)$ of the first $N$ terms of the sequence $\{S_k \alpha\}$, where $\{ \cdot \}$ denotes fractional part.
Under some additional regularity conditions (including $\E X_1 =0$ in the case $1< \beta <2$) we proved that
\begin{equation}\label{AA}
D_N (S_k \alpha ) = O \left( \sqrt{\frac{\log \log N}{N}} \log N \right), \quad D_N (S_k \alpha ) = \Omega \left( \sqrt{\frac{\log \log N}{N}} \right) \quad \text{a.s.}
\end{equation}
if $\beta \gamma \le 2$, and
 \begin{equation} \label{AA2} D_N (S_k \alpha ) = O \left( \left( \frac{\log \log N}{N} \right)^{1/(\beta \gamma)} \right), \quad D_N (S_k \alpha ) = \Omega \left( \frac{1}{N^{1/(\beta \gamma)}} \right) \quad \text{a.s.}
 \end{equation}
if $\beta \gamma >2$. Here a.s. (almost surely) means that a given relation holds with probability $1$. These relations show the surprising fact that the asymptotic behavior of $D_N(S_k\alpha)$ changes substantially at the critical value $\gamma= 2/\beta$. A similar result holds in the case $\E X_1^2<\infty$, $\E X_1=0$ (which corresponds formally to $\beta=2$)  when the critical value is $\gamma=1$, and in the case $\E |X_1| <\infty$, $\E X_1 \ne 0$, when the critical value is $\gamma=2$. (In the latter case the result was proved only for some special distributions.)
Note that there is a gap between the upper and lower bounds in (\ref{AA}) and (\ref{AA2}), and one purpose of the present paper is to remove the factor $\log N$ in the upper bound of (\ref{AA}), yielding a precise law of the iterated logarithm (LIL) for $D_N (S_k \alpha)$ in the case $\beta\gamma<1$. Under the same condition we will also prove that $N^{-1/2} D_N (S_k \alpha)$ has a nondegenerate limit distribution. In fact, we will prove the functional version of these results describing the asymptotic behavior of a large class of functionals of the sequence $\{S_k \alpha\}$, including the $L^p$-discrepancy $D_N^{(p)} (S_k \alpha )$, $1 \le p \le \infty$. Recall that $D_N^{(p)}(S_k \alpha)$ is defined as the $L^p ([0,1])$-norm of $F_N(t)-t$, where
\[ F_N (t) = \frac{1}{N} \sum_{k=1}^N I_{[0,t)} (\{ S_k \alpha \} ), \qquad 0 \le t \le 1 \]
is the empirical distribution function, and that the discrepancy can also be expressed in terms of $F_N$ as $D_N (S_k \alpha ) = \sup_{0 \le s<t \le 1} |(F_N(t)-t)-(F_N(s)-s)|$. We define $F_N$, $D_N^{(p)} (a_k)$ and $D_N (a_k)$ for an arbitrary sequence of reals $a_k$ analogously.

From a purely probabilistic point of view $\{S_k \alpha\}$ is a Markov chain with countable state space, but we can also view it on the state space $[0,1)$; in the latter case its stationary distribution is the uniform distribution on $[0,1)$. According to a classical result of L\'evy \cite{L}, the distribution of $\{S_k \alpha\}$ converges weakly to the uniform distribution as $k\to \infty$ if and only if $X_1$ is nondegenerate (i.e.\ attains more than one integers with positive probability). Since the limit distribution is continuous, we then have
\[ \psi_\alpha(k):= \sup_{0\le t \le 1} |\P (\{S_k \alpha\} \le t)-t| \longrightarrow 0 \ \, \text{as} \ \, k\to\infty . \]
Note that $D_N^{(\infty)}(S_k \alpha)$ is the empirical analogue of $\psi_{\alpha}(k)$. The order of magnitude of $\psi_{\alpha}(k)$ was investigated by Schatte \cite{SCH2}, Diaconis \cite{DI}, Su \cite{SU}, Hensley and Su \cite{HS}; improving their results, in \cite{BB} we showed that
\[ \psi_\alpha (k) =O(k^{-1/(\beta \gamma)}) . \]
Moreover, under mild additional regularity conditions on the distribution of $X_1$ this estimate is sharp. This shows that, in contrast to $D_N(S_k \alpha)$, there is no change of behavior of $\psi_{\alpha}(k)$ at $\beta \gamma=2$. Sharp results for the remainder term in the central limit theorem (CLT) and the corresponding Edgeworth expansion for i.i.d.\ sums (without taking the fractional part) under Diophantine conditions were given by Bobkov \cite{BOB1}, \cite{BOB2}.

By the previous result, condition $\beta \gamma<1$ amounts to the same as
$$ \psi_\alpha (k) \ll k^{-(1+\delta)} \quad  (\delta>0),$$
and in this paper we will extend our results on $D_N(S_k\alpha)$ to the discrepancy $D_N (S_k)$ with any i.i.d.\ sequence $X_1, X_2, \dots$ satisfying $\psi (k) \ll k^{-(1+\delta)}$ with some $\delta >0$, where
\begin{equation}\label{psi}
\psi (k): = \sup_{t \in [0,1]} |\P (\{ S_k\} \le t) -t| .
\end{equation}
As in the arithmetic case, we will also prove the functional version of our results, describing the asymptotics of general functionals of the sequence $\{S_k\}$. The LIL for $D_N(S_k)$ was proved by Schatte \cite{SCH3} under the assumption $\psi(k)\ll k^{-(5/2+\delta)}$, $\delta>0$; he also claimed, without proof, that the result remains valid under $\psi(k)\ll k^{-(2+\delta)}$, $\delta>0$. We believe that our condition $\psi (k) \ll k^{-(1+\delta)}$, $\delta >0$ is sharp, but its optimality remains open. For earlier results on the discrepancy of $D_N(S_k)$ in the case of absolutely continuous  $X_1$ (when $\psi(k)\to 0$ exponentially fast), see Bazarova et al. \cite{BBR}, Berkes and Raseta \cite{BR}.

To formulate our results, let $X_1, X_2, \dots$ be i.i.d.\ real-valued random variables, let $S_k=\sum_{j=1}^k X_j$, and let $\psi (k)$ be as in \eqref{psi}. Let $\mathcal{F}$ denote the set of all $1$-periodic functions $f: \R \to \R$ such that $f$ is of bounded variation on $[0,1]$ and $\int_0^1 f(x) \, \mathrm{d}x=0$. For any $f,g \in \mathcal{F}$ set
\begin{equation}\label{C}
C(f,g) = \E f(U) g(U) + \sum_{k=1}^{\infty} \E f(U)g(U+S_k) + \sum_{k=1}^{\infty} \E g(U)f(U+S_k),
\end{equation}
where $U$ is a random variable uniformly distributed on $[0,1)$, independent of $X_1, X_2, \dots$.
\begin{thm}\label{theorem1}
Assume that $\psi (k) \ll k^{-(1+\delta)}$ for some $\delta >0$. Then for any $f,g \in \mathcal{F}$ the series in \eqref{C} are absolutely convergent, and $C(f,g)$ is bilinear, symmetric and positive semidefinite in $f,g$; moreover, $C(f,f)=0$ if and only if $f=0$ a.e. For any $f \in \mathcal{F}$ the sum $\sum_{k=1}^N f(S_k)$ satisfies the central limit theorem
\begin{equation}\label{CLT}
\frac{\sum_{k=1}^N f(S_k)}{\sqrt{N}} \overset{d}{\longrightarrow} \mathcal{N} \left( 0,\sigma^2 \right)
\end{equation}
and the law of the iterated logarithm
\[ \limsup_{N \to \infty} \frac{\sum_{k=1}^N f(S_k)}{\sqrt{2N \log \log N}} =\sigma \quad \mathrm{a.s.} \]
with $\sigma = \sqrt{C(f,f)}$.
\end{thm}
\noindent Relation \eqref{CLT} means convergence in distribution to the mean zero normal distribution of variance $\sigma^2$ if $\sigma >0$, and to the constant zero if $\sigma =0$.

Next, let
\[ F_N (t) = \frac{1}{N} \sum_{k=1}^N I_{[0,t)} (\{ S_k \} ), \qquad 0 \le t \le 1 \]
denote the empirical distribution function of the sequence $\{S_k\}$. Further, for any interval $J \subseteq [0,1)$ of length $\lambda (J)$ let $f_J(x)=I_J(\{ x \})-\lambda (J)$ denote the centered indicator of $J$ extended with period $1$, and put $\Gamma (s,t) = C(f_{[0,s)}, f_{[0,t)})$, $s,t \in [0,1]$. The following two results describe the asymptotic behavior of $F_N(t)$ in the distributional and a.s.\ sense. For basic facts about functional limit theorems we refer to Billingsley \cite{BIL}.
\begin{thm}\label{theorem2} Assume that $\psi (k) \ll k^{-(1+\delta)}$ for some $\delta >0$.
Then the function $\Gamma (s,t)$ is symmetric, positive semidefinite and continuous on the unit square, and
\begin{equation}\label{donsker}
\sqrt{N} (F_N(t) - t) \overset{D[0,1]}{\longrightarrow} K(t) \quad  \text{as} \quad N \to \infty ,
\end{equation}
where $K(t)$ is a mean zero Gaussian process on $[0,1]$ with covariance function $\Gamma (s,t)$.
\end{thm}
\noindent Relation (\ref{donsker}) expresses weak convergence in the Skorohod space $D[0,1]$, and implies that the distribution of any functional $\Psi$ of the process $\sqrt{N} (F_N(t) - t)$, $0\le t \le 1$ converges weakly to the distribution of $\Psi (K(t))$, provided that $\Psi$ is continuous in the $D[0, 1]$ topology. This gives the limit distribution of a large class of functionals of the sequence $\{S_k\}$. Since the discrepancy and the $L^p$-discrepancy are continuous functionals, we immediately obtain
\[ \begin{split} \sqrt{N} D_N(S_k) &\overset{d}{\longrightarrow} \sup_{0\le s<t \le 1} |K(t)-K(s)| , \\ \sqrt{N} D_N^{(p)}(S_k) &\overset{d}{\longrightarrow} \| K \|_{L^p([0,1])} \qquad (1 \le p \le \infty ) . \end{split} \]
Note that $\Gamma (0,0)=\Gamma (1,1)=0$. From Theorem \ref{theorem1} it follows easily that $\Gamma$ is (strictly) positive definite on the open interval $(0,1)$; in particular, the limit distributions in the previous formula are nondegenerate.

The next result is the LIL analogue of Theorem \ref{theorem2}.
\begin{thm}\label{theorem3}
Assume that $\psi (k) \ll k^{-(1+\delta)}$ for some $\delta >0$. Then with probability $1$, the sequence of functions
\[ \frac{N (F_N(t)-t)}{\sqrt{2 N \log \log N}}, \quad 0\le t \le 1 \]
is relatively compact in the Skorohod space $D[0,1]$, and its class of limit functions in the $D[0,1]$ topology is identical with the closed unit ball $B(\Gamma )$ of the reproducing kernel Hilbert space determined by the covariance function $\Gamma (s,t)$.
\end{thm}
\noindent As a consequence, with probability $1$ we have
\[ \begin{split} \limsup_{N\to\infty} \frac{N D_N(S_k)}{\sqrt{2 N \log \log N}} &= \sup_{y \in B(\Gamma )} \sup_{0 \le s<t\le 1} |y(t)-y(s)|, \\ \limsup_{N\to\infty} \frac{N D_N^{(p)}(S_k)}{\sqrt{2 N \log \log N}} &= \sup_{y\in B(\Gamma )} \|y\|_{L^p([0,1])} \qquad (1 \le p \le \infty ) . \end{split} \]
Note that the value of these limsups are positive and finite. Indeed, finiteness is a consequence of the fact that $B(\Gamma)$ is compact in the supremum norm (see the proof of Theorem \ref{theorem3}), whereas positivity follows from Theorem \ref{theorem1} with $f=f_{[0,t)}$ for some $0<t<1$. Again, Theorem \ref{theorem3} yields laws of the iterated logarithm for a large class of functionals of the process $\sqrt{N} (F_N(t)-t)$, $0\le t \le 1$; see Strassen \cite{STR}.

The connection of reproducing kernel Hilbert spaces with the functional LIL for Gaussian processes was first observed by Oodaira \cite{OO}; analogous results for mixing processes and lacunary series were obtained by Philipp \cite{PH}. The latter paper provides a good introduction to the probabilistic aspects of reproducing kernel Hilbert spaces for the reader not familiar with this concept. The same paper also provides the background for the chaining techniques required in the proof of our theorems. For the general theory of reproducing kernel Hilbert spaces we refer to Aronszajn \cite{A}.

If $X_1$ is uniformly distributed on $[0,1)$, then it is not difficult to see that $\{ S_k \}$ is an i.i.d.\ sequence of uniformly distributed random variables. In this simplest case we have $C(f,g)=\E f(U)g(U)$, and consequently $\Gamma (s,t) = \min \{ s,t \} -st$ is the covariance function of the Brownian bridge. The corresponding functional CLT and functional LIL are due to Donsker \cite{DO} and Finkelstein \cite{F}.

In general, $\{S_k \}$ is a sequence of weakly dependent random variables, and there exist in the literature functional CLT's and LIL's for the empirical distribution function for such sequences, but their assumptions are considerably more restrictive than those for the central limit theorem. For example, for $\varphi$-mixing sequences with mixing rate $\varphi_n$, a functional CLT for the empirical distribution function holds under $\sum_{n=1}^\infty n^2 \varphi_n^{1/2}<\infty$  (see Billingsley \cite[p.\ 197, Theorem 22.1]{BIL}), while the (ordinary) CLT for partial sums holds under $\sum_{n=1}^\infty \varphi_n^{1/2}<\infty$. Also, most results in the literature deals with stationary processes, while $\{S_k \}$ is nonstationary.

\section{Couplings and the Kolmogorov metric}

We start with two useful observations on the relationship between the Kolmogorov metric, functions of bounded variation and couplings. The total variation of a function $f$ on $[0,1]$ is denoted by $V(f)$, and the $L^p([0,1])$-norm by $\| f \|_p$.
\begin{lem}\label{randomkoksma}
For any $[0,1]$-valued random variables $X$ and $Y$, and any function $f: [0,1] \to \mathbb{R}$ of bounded variation,
\[ \left| \E f(X)- \E f(Y) \right| \le V(f) \sup_{0 \le t \le 1} \left| \P ( X \le t) -\P (Y \le t) \right| . \]
\end{lem}
\noindent In the special case when the distribution of $X$ is finitely supported with equal weights and $Y$ is uniformly distributed on $[0,1]$, this is known as Koksma's inequality \cite[Chapter 2, Theorem 5.1]{KN}. In particular, for any real-valued random variables $X$ and $Y$, and any $f \in \mathcal{F}$ we have
\begin{equation}\label{maintool}
\left| \E f (X) - \E f (Y) \right| \le V(f) \sup_{0 \le t \le 1} \left| \P ( \{ X \} \le t) -\P (\{Y \} \le t) \right| .
\end{equation}
The claim formally follows from integration by parts; however, as one-sided continuity of $f$ is not assumed, integration by parts strictly speaking cannot be applied.
\begin{proof}[Proof of Lemma \ref{randomkoksma}] First, let $g:[0,1] \to \R$ be a nondecreasing function with $g(0)=0$. Note that for any $x \in [0,g(1)]$ the set $L(x)=\{ y \in [0,1] : g(y) > x \}$ is an interval of the form $[t,1]$ or $(t,1]$ for some $0 \le t \le 1$. Therefore
\[ \begin{split} \left| \E g(X) -\E g(Y) \right| &= \left| \int_0^{g(1)} \left( \P (g(X) >x) - \P (g(Y) >x) \right) \, \mathrm{d}x \right| \\ &= \left| \int_0^{g(1)} \left( \P ( X \in L(x)) - \P (Y \in L(x)) \right) \, \mathrm{d}x \right| \\ &\le g(1) \sup_{0 \le t \le 1} \left| \P (X \le t) -\P (Y \le t) \right| . \end{split} \]
Here $g(1)=V(g)$. Since the left hand side and $V(g)$ are both invariant under adding a constant to $g$, we have
\begin{equation}\label{nondecreasing}
\left| \E g(X) -\E g(Y) \right| \le V(g) \sup_{0 \le t \le 1} \left| \P (X \le t) -\P (Y \le t) \right|
\end{equation}
for any nondecreasing function $g:[0,1] \to \R$. Finally, consider the Jordan decomposition $f=g_1-g_2$, where $g_1, g_2$ are nondecreasing functions and $V(f)=V(g_1)+V(g_2)$. Applying \eqref{nondecreasing} to $g_1$ and $g_2$, the claim follows.
\end{proof}

\begin{lem}\label{couplinglemma} Let $X$ and $\xi$ be independent $[0,1]$-valued random variables, and assume that $\xi$ is uniformly distributed on $[0,1]$. Then there exists a uniformly distributed $[0,1]$-valued random variable $U$ such that $U$ is a function of $X$ and $\xi$, and $|U-X| \le \sup_{0 \le t \le 1} |\P (X \le t) -t|$.
\end{lem}

\begin{proof} Let $F(x)=\P (X\le x)$ and $F^-(x)=\P (X < x)$, and define
\[ U=\xi F(X)+(1-\xi )F^-(X). \]
We have $\left| U-X \right| \le \sup_{0 \le t \le 1} \left| F(t)-t \right|$ by construction. To show that $U$ is uniformly distributed on $[0,1]$, let $D$ denote the set of points at which $F$ is not continuous. If $x \in [0,1] \backslash \bigcup_{d \in D} [F^{-}(d),F(d))$, then $t= \sup \{ u \in [0,1] : F(u) \le x \}$ is seen to satisfy $x=F(t)$. Thus $\{ U \le x \} = \{ F(X) \le F(t) \} = \{ X \le t \}$ as random events, and so $\P (U \le x ) = \P (X \le t) = F(t)=x$. If $x \in [F^{-}(d),F(d))$ for some $d \in D$, then
\[ \left\{ U \le x \right\} = \{ X < d \} \cup \left\{ X=d \textrm{ and } \xi \le \frac{x-F^{-}(d)}{F(d)-F^{-}(d)} \right\} \]
as random events, hence using the independence of $X$ and $\xi$ we get
\[ \begin{split} \P (U \le x) &= \P (X<d) + \P (X=d) \P \left( \xi \le \frac{x-F^{-}(d)}{F(d)-F^{-}(d)} \right) \\ &=F^{-}(d) + \left( F(d)-F^{-}(d) \right) \cdot \frac{x-F^{-}(d)}{F(d)-F^{-}(d)} =x. \end{split} \]
We note that if $X$ has a continuous distribution, then $F=F^{-}$ and so $U=F(X)$ has the same property; a classical observation of probability theory.
\end{proof}

Given two $[0,1]$-valued random variables $X$ and $Y$ with distribution functions $F(t)=\P (X \le t)$ and $G(t)=\P (Y \le t)$, their distance in the Kolmogorov metric (or uniform metric) is defined as
\[ d_{\mathrm{Kol}}(X,Y) = \sup_{0 \le t \le 1} |\E I_{[0,t]}(X) - \E I_{[0,t]}(Y)| = \sup_{0 \le t \le 1}|F(t)-G(t)|. \]
Lemma \ref{randomkoksma} states that here the family of test functions $I_{[0,t]}$, $0 \le t \le 1$ can be replaced by the larger class of functions of bounded variation; that is,
\[ d_{\mathrm{Kol}}(X,Y) = \sup_{V(f) \le 1} |\E f(X) - \E f(Y)| . \]
Another classical distance for $[0,1]$-valued random variables is the $p$-Wasserstein metric, defined as
\[ W_p (X,Y) = \inf \left\{ \| X'-Y' \|_p \, : \, X' \overset{d}{=} X \textrm{ and } Y' \overset{d}{=}Y \right\} \qquad (1 \le p \le \infty ). \]
Here $X' \overset{d}{=} X$ means that $X'$ and $X$ have the same distribution; the infimum is thus over all couplings $(X', Y')$ of $X$ and $Y$. The $p$-Wasserstein metric can be expressed in terms of the distribution functions \cite[Theorem 2.10]{BL} as $W_p (X,Y) = \| F^{-1} - G^{-1} \|_p$, $1 \le p \le \infty$, where $F^{-1} (x) = \inf \{ t \in [0,1] \, : \, F(t) \ge x \}$, $0<x<1$ is the generalized inverse of $F$.

It is not difficult to see that for $G(t)=t$,
\[ \| F^{-1} -G^{-1} \|_{\infty} = \sup_{0<x<1} |F^{-1}(x) -x| = \sup_{0<t<1} |F(t) -t|. \]
Indeed, by the left-continuity of $F^{-1}$ the essential supremum and the actual supremum on $0<x<1$ are equal; to see the second equality note that $|F(t)-t| \le c$ geometrically means that the graph of $F(t)$, $0<t<1$ stays between the lines $t \pm c$, and that this is invariant under taking the inverse function. Since in the definition of the Kolmogorov metric the supremum can also be taken over $0<t<1$, it follows that $W_{\infty} (X,U) = d_{\mathrm{Kol}}(X,U)$ for a uniformly distributed $U$; Lemma \ref{couplinglemma} is basically a constructive form of this fact. We summarize these observations in the following lemma.
\begin{lem} Let $X$ and $U$ be $[0,1]$-valued random variables, and assume that $U$ is uniformly distributed. Then
\[ W_{\infty} (X,U) = \sup_{V(f) \le 1} |\E f(X) - \E f(U)|, \]
where the supremum is over all functions of bounded variation $f: [0,1] \to \mathbb{R}$ of total variation $V(f) \le 1$.
\end{lem}
\noindent This duality between $W_{\infty}$ and functions of bounded variation seems not to have been observed before. The classical Kantorovich duality theorem \cite[Theorem 2.5]{BL}
\[ W_1 (X,Y) = \sup_{\mathrm{Lip} (f) \le 1} |\E f(X) - \E f(Y) | , \]
where the supremum is over all $f:[0,1] \to \mathbb{R}$ such that $|f(x)-f(y)| \le |x-y|$, establishes a duality between $W_1$ and Lipschitz functions. Amusingly, the general inequality $W_1 (X,U) \le W_{\infty} (X,U)$ thus reduces to the fact that every Lipschitz function is of bounded variation.

\section{The variance}

We now find the precise asymptotics of the variance.
\begin{lem}\label{varianceprop1} Suppose that $\psi (k) \ll k^{-(1+\delta)}$ for some $\delta >0$.
\begin{enumerate}
\item[(i)] For any $f,g \in \mathcal{F}$ the series in \eqref{C} are absolutely convergent, and $|C(f,g)| \ll V(f) \| g \|_1 + V(g) \| f \|_1$. Further, $C(f,g)$ is bilinear, symmetric and positive semidefinite in $f,g$. In addition, $C(f,f)=0$ if and only if $f=0$ a.e.
\item[(ii)] For any $f \in \mathcal{F}$ and any integers $M \ge 0$ and $N \ge 1$,
\[ \begin{split} \E \left( \sum_{k=M+1}^{M+N} f(S_k) \right)^2 = C(f,f) N +O \Big( &V(f) \left( \| f \|_1 + V(f) (M+1)^{-\delta} \right)^{\delta /(1+\delta )} \times \\ &\left( \| f \|_1 N + V(f) (M+1)^{-\delta} \right)^{1/(1+\delta )} \Big) . \end{split} \]
\end{enumerate}
The implied constants in (i) and (ii) depend only on the distribution of $X_1$ and $\delta$.
\end{lem}

\begin{proof} First, we prove (i). In the definition \eqref{C} of $C(f,g)$ the variable $U$ is independent of $X_1, X_2, \dots$, therefore we can write
\[ \E f(U)g(U+S_k) = \E \int_0^1 f(u)g(u+S_k) \, \mathrm{d}u = \E G(S_k), \]
where $G(x)=\int_0^1 f(u)g(u+x) \, \mathrm{d} u$. It is not difficult to see e.g.\ directly from the definition of total variation that $V(G)\le V(g) \| f \|_1$. Applying \eqref{maintool} to $X=S_k$ and a random variable $Y$ uniformly distributed on $[0,1)$ and noting that $\E G(Y)=0$, we get $|\E f(U)g(U+S_k)| \le V(g) \| f \|_1 \psi (k)$. Repeating the same argument for $\E g(U) f(U+S_k)$, the absolute convergence of the series in \eqref{C} follows.  Since $|\E f(U)g(U)| \le V(g) \| f \|_1$, we also get $|C(f,g)|\ll V(f) \| g \|_1 + V(g) \| f \|_1$. Clearly $C(f,g)$ is bilinear and symmetric in $f,g$.

Let $\hat{f}(h)=\int_0^1 f(x) e^{-2 \pi i h x} \, \mathrm{d}x$, $h \in \mathbb{Z}$ denote the Fourier coefficients of $f$. We claim that
\begin{equation}\label{fouriervariance}
C(f,f)=\sum_{h \neq 0} |\hat{f}(h)|^2 \frac{1-|\varphi (2 \pi h)|^2}{|1-\varphi (2 \pi h)|^2} ,
\end{equation}
where $\varphi (x)=\E \exp (i x X_1)$ is the characteristic function of $X_1$. Here $|\varphi (2 \pi h)|<1$ for all integers $h \neq 0$; indeed, otherwise $X_1$ would have a lattice distribution with a rational span, which is impossible by the assumption $\psi (k) \to 0$. From \eqref{fouriervariance} it will follow immediately that $C(f,f) \ge 0$, and that $C(f,f)=0$ if and only if $\hat{f}(h)=0$ for all integers $h \neq 0$; the latter condition is equivalent to $f=0$ a.e. We will also have
\[ \sum_{i,j=1}^r C(f_i,f_j) x_i x_j = C \left( \sum_{i=1}^r f_i x_i, \sum_{i=1}^r f_i x_i \right) \ge 0 \]
for any $f_1, f_2, \dots, f_r \in \mathcal{F}$ and any $x_1, x_2, \dots, x_r \in \mathbb{R}$, proving that $C(f,g)$ is positive semidefinite.

To see \eqref{fouriervariance}, assume first that $f$ is a trigonometric polynomial; that is, $\hat{f}(h)=0$ for all but finitely many integers $h$. As before, $\E f(U)f(U+S_k)=\E G(S_k)$ with $G(x)=\int_0^1 f(u)f(u+x) \, \mathrm{d}u$. Since $\hat{G}(h)=|\hat{f}(h)|^2$, we have $G(x)=\sum_{h \neq 0} |\hat{f}(h)|^2 e^{2 \pi i hx}$. Therefore
\[ \sum_{k=1}^{\infty} \E G(S_k) = \sum_{h \neq 0} \sum_{k=1}^{\infty} |\hat{f}(h)|^2 \varphi (2 \pi h)^k = \sum_{h \neq 0} |\hat{f}(h)|^2 \frac{\varphi (2 \pi h)}{1-\varphi (2 \pi h)} .  \]
Combining the $h$ and $-h$ terms, and using $\E f(U)^2 = \sum_{h \neq 0} |\hat{f}(h)|^2$, we obtain \eqref{fouriervariance}. For a general $f \in \mathcal{F}$, consider
\[ f_H (x) = \int_0^1 f(x-y) F_H(y) \, \mathrm{d}y = \sum_{0<|h|<H} \left( 1-\frac{|h|}{H} \right) \hat{f}(h) e^{2 \pi i hx} , \]
where $F_H(x)=\sum_{|h|<H}(1-|h|/H)e^{2 \pi i hx}$ is the Fej\'er kernel of order $H$. Note that $f_H$ is a trigonometric polynomial; in fact, it is a Ces\`aro mean of the formal Fourier series of $f$. By the special case proved above, we have
\[ C(f_H, f_H) = \sum_{0<|h|<H} \left( 1-\frac{|h|}{H} \right)^2 |\hat{f}(h)|^2 \frac{1-|\varphi (2 \pi h)|^2}{|1-\varphi (2 \pi h)|^2} \to \sum_{h \neq 0} |\hat{f}(h)|^2 \frac{1-|\varphi (2 \pi h)|^2}{|1-\varphi (2 \pi h)|^2} \]
as $H \to \infty$. On the other hand, letting $G(x)=\int_0^1 f(u)f(u+x) \, \mathrm{d}u$ and $G_H(x)=\int_0^1 f_H(u)f_H(u+x) \, \mathrm{d}u$, by Lemma \ref{randomkoksma} we have
\[ \begin{split} \left| C(f,f) - C(f_H, f_H) \right| &= \left| \E f(U)^2 + 2 \sum_{k=1}^{\infty} \E G(S_k) - \E f_H(U)^2 - 2 \sum_{k=1}^{\infty} \E G_H (S_k) \right| \\ &\le \left| \E f(U)^2 - \E f_H(U)^2 \right| + 2 \sum_{k=1}^{\infty} |\E (G-G_H)(S_k)| \\ &\le \left| \| f \|_2^2 - \| f_H \|_2^2 \right| + 2 \sum_{k=1}^{\infty} V(G-G_H) \psi (k) \\ &\ll \left| \| f \|_2^2 - \| f_H \|_2^2 \right| + V(G-G_H) . \end{split} \]
Since $f_H \to f$ in $L^2([0,1])$, it will be enough to prove $V(G-G_H) \to 0$ as $H \to \infty$. Writing $f=f_H+(f-f_H)$ in the definition of $G(x)$ and using the integral transformation $u \mapsto u-x \pmod{1}$ in one of the terms, we get
\[ (G-G_H)(x) = \int_0^1 (f-f_H)(u) \left( f_H(u-x) + f(u+x) \right) \, \mathrm{d} u , \]
and consequently $V(G-G_H) \le \| f-f_H \|_1 \left( V(f_H) + V(f) \right)$. Using $F_H \ge 0$ it is readily seen that $V(f_H) \le V(f)$. Hence $V(G-G_H) \le 2V(f) \| f-f_H \|_1 \to 0$, finishing the proof of \eqref{fouriervariance} for a general $f \in \mathcal{F}$.

Next, we prove (ii). We may assume $V(f)=1$. Let $U$ be as in the definition \eqref{C} of $C(f,f)$. Expanding the square in the claim we get
\begin{equation}\label{expsquare}
\E \left( \sum_{k=M+1}^{M+N} f(S_k) \right)^2 = \sum_{k=M+1}^{M+N} \E f(S_k)^2 + 2 \sum_{M+1 \le k< \ell \le M+N } \E f(S_k) f(S_{\ell}) .
\end{equation}
To estimate the diagonal terms, let us apply \eqref{maintool} to $X=S_k$ and $Y=U$. Noting that $V(f^2) \le 2 V(f)^2=2$, we obtain $\E f(S_k)^2 = \E f(U)^2+O(\psi (k))$. The assumption $\psi (k) \ll k^{-(1+\delta)}$ thus ensures that $\sum_{k=M+1}^{M+N} \E f(S_k)^2 = N \E f(U)^2 +O\left( (M+1)^{-\delta} \right)$.

We now estimate the off-diagonal terms in two different ways. For given $k<\ell$ let us write $f(S_k) f(S_{\ell}) = f(S_k) f(S_k +(S_{\ell}-S_k))$. First, let $\mathcal{F}_k$ denote the $\sigma$-algebra generated by $X_1, X_2, \dots, X_k$, and note that $S_{\ell}-S_k$ is independent of $\mathcal{F}_k$. Let us apply \eqref{maintool} to $X=S_{\ell}-S_k$, a random variable $Y$ uniformly distributed on $[0,1)$ independent of $\mathcal{F}_k$, and the function $g(x)=f(S_k) f(S_k+x)$ to estimate the conditional expectation with respect to $\mathcal{F}_k$. Since $\E (g(Y) | \mathcal{F}_k)=0$ and $V(g) \le |f(S_k)| V(f)=|f(S_k)|$, we have
\[ \left| \E \left( f(S_k) f(S_{\ell}) \mid \mathcal{F}_k \right) \right| \le |f(S_k)| \psi (\ell -k) . \]
Taking the (total) expectation, we obtain $|\E f(S_k) f(S_{\ell})| \le \E |f(S_k)| \psi (\ell -k)$. We can apply \eqref{maintool} again, this time to $X=S_k$, a random variable $Y$ uniformly distributed on $[0,1)$ and $|f|$. Since $\E |f(Y)|=\| f \|_1$ and $V(|f|) \le V(f) = 1$, we get
\begin{equation}\label{skslfirst}
|\E f(S_k) f(S_{\ell})| \le (\| f \|_1 + \psi (k)) \psi (\ell -k).
\end{equation}
Second, let $\mathcal{G}_k$ denote the $\sigma$-algebra generated by $X_{k+1}, X_{k+2}, \dots$. Note that $S_k$ and $U$ are independent of $\mathcal{G}_k$. Let us apply \eqref{maintool} to $X=S_k$, $Y=U$ and the function $g(x)=f(x)f(x+(S_{\ell}-S_k))$. Since $V(g) \le 2 V(f)^2=2$, we obtain
\[ \E \left( f(S_k) f(S_{\ell}) \mid \mathcal{G}_k \right) = \E \left( f(U)f(U+(S_{\ell}-S_k)) \mid \mathcal{G}_k \right) + O(\psi (k)) . \]
Note that $f(U)f(U+(S_{\ell}-S_k))$ has the same distribution as $f(U)f(U+S_{\ell -k})$. Taking the (total) expectation, we thus get
\begin{equation}\label{skslsecond}
\E f(S_k) f(S_{\ell}) = \E f(U) f(U+S_{\ell -k})+O(\psi (k)).
\end{equation}
As observed in (i), here $|\E f(U) f(U+S_{\ell -k})| \le \| f \|_1 \psi (\ell -k)$.

Let $1 \le K \le N$ be an integer, to be chosen. We will estimate the off-diagonal terms for which $\ell -k > K$, and those for which $\ell -k \le K$ separately. The assumption $\psi (k) \ll k^{-(1+\delta)}$ and \eqref{skslfirst} show that
\[ \begin{split} \sum_{\substack{M+1 \le k< \ell \le M+N \\ \ell -k >K}} \E f(S_k) f(S_{\ell}) &\ll \left( \| f \|_1 N+(M+1)^{-\delta} \right) K^{-\delta}, \\ \sum_{M+N-K<k< \ell \le M+N} \E f(S_k) f(S_{\ell}) &\ll \| f \|_1 K + (M+1)^{-\delta}. \end{split} \]
For a fixed $M+1 \le k \le M+N-K$, from \eqref{skslsecond} we obtain
\[ \begin{split} \sum_{\ell =k+1}^{k+K} \E f(S_k) f(S_{\ell}) &= \sum_{d=1}^K \E f(U) f(U+S_d) + O \left( K \psi (k) \right) \\ &= \sum_{d=1}^{\infty} \E f(U) f(U+S_d) + O \left( \| f \|_1 \sum_{d=K+1}^{\infty} \psi (d) + K \psi (k) \right) . \end{split} \]
Summing over $M+1 \le k \le M+N-K$ and using the assumption $\psi (k) \ll k^{-(1+\delta)}$, we altogether get
\begin{multline*}
\sum_{M+1 \le k< \ell \le M+N} \E f(S_k) f(S_{\ell}) = \\ N \sum_{d=1}^{\infty} \E f(U) f(U+S_d)+O \left( \left( \| f \|_1 N +(M+1)^{-\delta} \right) K^{-\delta} + \left( \| f \|_1 + (M+1)^{-\delta} \right) K \right) .
\end{multline*}
Choosing $K = \lfloor \left( \| f \|_1 N +(M+1)^{-\delta} \right)^{1/(1+\delta)} / \left( \| f \|_1 +(M+1)^{-\delta} \right)^{1/(1+\delta)} \rfloor$, in \eqref{expsquare} we thus have
\begin{multline*}
\E \left( \sum_{k=M+1}^{M+N} f(S_k) \right)^2 = \\ C(f,f) N+ O\left( \left( \| f \|_1 + (M+1)^{-\delta} \right)^{\delta /(1+ \delta )} \left( \| f \|_1 N + (M+1)^{-\delta} \right)^{1/(1+\delta )} \right) .
\end{multline*}
\end{proof}

\section{Approximation by independent variables}\label{approximationsubsection}

In this section we approximate $\sum_{k=M+1}^{M+N} f(S_k)$ by a sum of independent random variables. A basic tool is the following result, implicit in Schatte \cite{SCH3}, \cite{SCH4}.
\begin{lem}\label{schattelemma}
Let $\xi_1, \xi_2, \ldots$ be independent random variables uniformly distributed on $[0,1)$, independent of the $X_n$'s, and let $I_1, I_2, \ldots$ be disjoint blocks of integers such that between $I_n$ and $I_{n+1}$ there are $\ell_n\ge 1$ integers. Then there exists a sequence $\delta_1, \delta_2, \dots$ of random variables such that $\delta_n$ is measurable with respect to $\xi_n$ and the $X_j$'s between $I_n$ and $I_{n+1}$, and $|\delta_n| \leq \psi(\ell_n)$, $n=1, 2, \ldots$; further, the random vectors
$$(\{S_i\}, i\in I_1), \  (\{S_i - \delta_1 \}, i \in I_2),  \ldots,  (\{S_i - \delta_{n - 1}\}, i \in I_n), \ldots $$
are independent and they have, except for the first one, uniformly distributed components.
\end{lem}

\begin{proof} We begin with a simple remark. Let $X$ be a random variable and let $U$ be a uniformly distributed random variable independent of $X$. Then $\{ X+U \}$ is uniformly distributed on $[0,1)$ and is independent of $X$, as one can see immediately by conditioning $\{ X+U \}$ on $\{X=c\}$.

Let now
$$ I_n := \left\{ p_n +1, p_n + 2, \dots, p_n + b_n\right\}, \quad n=1, 2, \ldots $$
We claim that
\[ \delta_n =\{ S_{p_{n + 1}} - S_{p_n + b_n} \} -U (\{ S_{p_{n + 1}}  - S_{p_n + b_n}\} , \xi_n ), \quad n=1, 2, \ldots \]
satisfies the lemma, where $U (\{ S_{p_{n + 1}} - S_{p_n + b_n} \}, \xi_n )$ is the uniform random variable constructed in Lemma \ref{couplinglemma} for the variables $\{ S_{p_{n + 1}} - S_{p_n + b_n} \}$ and $\xi_n$. Since $S_{p_{n+1}} - S_{p_n+ b_n} \stackrel{d}{=} S_{\ell_n}$, by the definition of $\psi (k)$ and the claim of Lemma \ref{couplinglemma}, $| \delta_n| \leq \psi(\ell_n)$.

Furthermore, we have
\begin{align*}
\{ S_{p_2 + 1} - \delta_1 \} &=
\{ \left(X_1 + \dots + X_{p_1+b_1}\right) +X_{p_2+1} +U (\{ S_{p_2} -S_{p_1+b_1}\}, \xi_1 ) \} \\
\vdots & \\
\{ S_{p_2 + b_2} - \delta_1 \} & =\{ \left( X_1 + \dots + X_{p_1+b_1}\right) +(X_{p_2+1}+ \cdots + X_{p_2+b_2}) \\
\phantom{9999999} & + U (\{ S_{p_2} -S_{p_1+b_1}\} , \xi_1 ) \} .
\end{align*}
Put
\begin{align*}
V &= \{ X_1 + \dots + X_{p_1+b_1}+ U (\{ S_{p_2} - S_{p_1+b_1}\}, \xi_1 )\},\\
X &= (X_{p_2+1}, X_{p_2+1} +X_{p_2+2}, \ldots, X_{p_2+1} + \ldots +X_{p_2+b_2} ).
\end{align*}
Let $\omega$ be a point of the underlying probability space, and put
\[ c=(X_1 + \dots + X_{p_1+b_1})(\omega ). \]
The conditional distribution of the vector $(\{S_{p_2+1}-\delta_1\}, \ldots, \{S_{p_2+b_2}-\delta_1\})$ relative to $\sigma\{X_1, \ldots, X_{p_1+b_1}\}$, evaluated at $\omega$, equals the distribution of
\[ \{X+c+ U (\{S_{p_2} - S_{p_1+b_1}\}, \xi_1)\} ,\]
where adding the scalar $V^*=c+ U (\{S_{p_2} - S_{p_1+b_1}\}, \xi_1)$ to $X$ is meant coordinatewise. But $U (\{S_{p_2} - S_{p_1+b_1}\}, \xi_1)$ is a uniform random variable independent of $X$ and so is $\{V^*\}$, showing that the distribution of $\{X+V^*\}$ does not depend on $c$ and consequently, $(\{S_{p_2+1}-\delta_1\}, \ldots, \{S_{p_2+b_2}-\delta_1\})$ is independent of $\sigma\{X_1, \ldots, X_{p_1+b_1}\}$ and thus of $(\{S_i\}, \, i\in I_1)$. By our remark at the beginning of the proof, the coordinates of $(\{S_{p_2+1} -\delta_1\}, \ldots, \{S_{p_2+b_2} -\delta_1\})$ are uniformly distributed.

Now let $n \ge 2$. As before, the vector
\begin{equation}\label{vectorn}
(\{ S_{p_{n+1}+1} -\delta_n \}, \ldots, \{ S_{p_{n+1}+b_{n+1}} -\delta_n \} )
\end{equation} can be written as $\{ X+V \}$, where
\begin{align*}
V &= \{ X_1 + \dots + X_{p_n+b_n} + U (\{ S_{p_{n+1}} - S_{p_n+b_n} \} , \xi_n ) \} ,\\
X &= (X_{p_{n+1}+1}, X_{p_{n+1}+1} +X_{p_{n+1}+2} , \ldots, X_{p_{n+1}+1} + \cdots +X_{p_{n+1}+b_{n+1}} ).
\end{align*}
Repeating the argument used in the case $n=1$, it follows that the vector (\ref{vectorn}) is independent of $\sigma\{X_1, \ldots, X_{p_n+b_n}, \xi_1, \ldots, \xi_{n-1} \}$
and thus of
$$(\{S_i\}, i\in I_1), \  (\{S_i - \delta_1\}, i \in I_2),  \ldots,  (\{S_i - \delta_{n - 1}\}, i \in I_n),
$$
and has uniform coordinates. This completes the proof of the lemma.
\end{proof}

Fix a function $f \in \mathcal{F}$, an integer $M \ge 0$ and a sequence of auxiliary random variables $\xi_2, \xi_3, \ldots$ uniformly distributed on $[0,1)$, independent of $X_1, X_2, \ldots$. Consider a decomposition of the infinite set of integers $\{ M+1, M+2, \dots \}$ into nonempty, consecutive blocks $H_1, J_1, H_2, J_2, \ldots$, and let $T_i=T_i^{(f)}=\sum_{k \in H_i} f(S_k)$ and $D_i=D_i^{(f)}=\sum_{k \in J_i} f(S_k)$ for all $i \ge 1$. If $M+N=\max J_R$ for some $R \ge 1$, we thus have
\[ \sum_{k=M+1}^{M+N} f(S_k) = \sum_{i=1}^R T_i + \sum_{i=1}^R D_i . \]
Using Lemma \ref{schattelemma} we can introduce two sequences of random variables $\delta_2, \delta_3, \ldots$ and $\delta_2', \delta_3', \ldots$ with the following properties. First, $\delta_i$ is measurable with respect to the $\sigma$-algebra generated by $X_k$, $k \in J_{i-1}$ and $\xi_i$, and $|\delta_i| \le \psi (|J_{i-1}|)$ for all $i \ge 2$; similarly, $\delta_i'$ is measurable with respect to the $\sigma$-algebra generated by $X_k$, $k \in H_i$ and $\xi_i$, and $|\delta_i'| \le \psi (|H_i|)$ for all $i \ge 2$. Letting $T_i^*=T_i^{(f)*}=\sum_{k \in H_i} f(S_k-\delta_i)$ and $D_i^*=D_i^{(f)*}=\sum_{k \in J_i} f(S_k-\delta_i')$, we also have that $T_2^*, T_3^*, \ldots$ are independent and $\E T_i^*=0$ for all $i \ge 2$; similarly, $D_2^*, D_3^*, \ldots$ are independent and $\E D_i^*=0$ for all $i \ge 2$. Under the condition $\psi (k) \ll k^{-(1+\delta)}$ the variances are given by
\begin{equation}\label{ti*di*variance}
\begin{split} \E T_i^{*2}&=C(f,f) |H_i|+O\left( V(f)^2 |H_i|^{1/(1+\delta )} \right), \\ \E D_i^{*2}&= C(f,f) |J_i|+O\left( V(f)^2 |J_i|^{1/(1+\delta )} \right) \end{split}
\end{equation}
with implied constants depending only on the distribution of $X_1$ and $\delta$. Indeed, conditioning on the $\sigma$-algebra $\mathcal{A}_i$ generated by the variables $X_k$, $k \in [1,M] \cup H_1 \cup J_1 \cup \cdots \cup H_{i-1} \cup J_{i-1}$ and $\xi_i$, we can write $T_i^*=\sum_{\ell=1}^{|H_i|} f(\tilde{S_\ell}-Z_i)$ with some $\mathcal{A}_i$-measurable random variable $Z_i$, where $\tilde{S_{\ell}}$ are the partial sums of the i.i.d.\ sequence $X_k$, $k \in H_i$. Applying Lemma \ref{varianceprop1} to $\tilde{f}(x)=f(x-Z_i)$ and noting that $C(\tilde{f},\tilde{f})=C(f,f)$ and $V(\tilde{f})=V(f)$, we obtain $\E \left( T_i^{*2} \mid \mathcal{A}_i \right) = C(f,f) |H_i|+O( V(f)^2 |H_i|^{1/(1+\delta )})$. Taking the (total) expectation, the first relation in \eqref{ti*di*variance} follows; the proof of the second relation is analogous.

\begin{lem}\label{approxerror1} Assume that $\psi (k) \ll k^{-(1+\delta)}$ with some $0<\delta<1$. If we choose the sizes of the blocks as $|H_i|=\lceil i^c \rceil$ and $|J_i|=\lceil i^{c'} \rceil$ for all $i \ge 1$, where $1/2<c<1$ is a fixed constant and $c'=(c+\delta/2)/(1+\delta)$, then
\[ \E \left( \sum_{i=2}^R (T_i - T_i^*) \right)^2 \ll V(f)^2 R^{1-\delta /2}, \quad \E \left( \sum_{i=2}^R (D_i - D_i^*) \right)^2 \ll V(f)^2 R^{1-\delta /2} \]
for all $R \ge 2$ with implied constants depending only on the distribution of $X_1$, $\delta$ and $c$.
\end{lem}

\begin{proof} We shall only give a proof for $\sum_{i=2}^R (T_i-T_i^*)$, as the proof for $\sum_{i=2}^R (D_i-D_i^*)$ is analogous. We may assume $V(f)=1$. Note that $c'<c$. Expanding the square we get
\begin{equation}\label{expandsq}
\E \left( \sum_{i=2}^R (T_i - T_i^*) \right)^2 = \sum_{i=2}^R \E (T_i-T_i^*)^2 + 2 \sum_{2 \le i<j\le R} \E (T_i-T_i^*)(T_j-T_j^*) .
\end{equation}
Let us first estimate the diagonal term $\E (T_i-T_i^*)^2$. Let $\mathcal{F}_i$ denote the $\sigma$-algebra generated by the variables $X_k$, $k \in J_{i-1}$ and $\xi_i$, and let $Y_i=\sum_{k \in J_{i-1}} X_k$. Note that in $T_i-T_i^*=\sum_{k \in H_i} (f((S_k-Y_i)+Y_i)-f((S_k-Y_i)+Y_i-\delta_i))$ the variables $Y_i$ and $\delta_i$ are $\mathcal{F}_i$-measurable, and $(S_k-Y_i)$ is independent of $\mathcal{F}_i$. Hence we can apply Lemma \ref{varianceprop1} to the i.i.d.\ sequence obtained from $X_1, X_2, \dots$ by deleting the terms $X_k$, $k \in J_{i-1}$ to estimate the conditional expectation with respect to $\mathcal{F}_i$. It is not difficult to see e.g.\ directly from the definition of total variation that the function $g(x)=f(x+Y_i)-f(x+Y_i-\delta_i)$ satisfies $V(g) \le 2 V(f)=2$ and $\| g \|_1 \ll |\delta_i|$. Therefore
\begin{multline*} \E \left( (T_i-T_i^*)^2 \mid \mathcal{F}_i \right) \ll \\ |\delta_i | \cdot |H_i| + \left( |\delta_i | + (\min H_i - |J_{i-1}|)^{-\delta} \right)^{\delta / (1+\delta )} \left( |\delta_i | \cdot |H_i| + (\min H_i - |J_{i-1}|)^{-\delta} \right)^{1/(1+ \delta )}.
\end{multline*}
Here $\min H_i - |J_{i-1}| \gg i^{c+1}$, and $|\delta_i| \le \psi (|J_{i-1}|) \ll i^{-c-\delta /2}$. Taking the (total) expectation, after simplifying the exponent we get $\E (T_i-T_i^*)^2 \ll i^{-\delta /2}$.
Summing over $2 \le i \le R$ we get that the contribution of the diagonal terms in \eqref{expandsq} is $\sum_{i=2}^R \E (T_i-T_i^*)^2 \ll R^{1-\delta/2}$.

We now estimate the off-diagonal terms in \eqref{expandsq} similarly to the proof of Lemma  \ref{varianceprop1}. For given $i<j$ let us write
\[ \begin{split} (T_i-T_i^*)(T_j-T_j^*) = \sum_{\substack{k \in H_i \\ \ell \in H_j}} &\left( f(S_k)-f(S_k -\delta_i ) \right) \times \\ &\left( f((S_{\ell}-W_{i,j})+W_{i,j})-f((S_{\ell}-W_{i,j})+W_{i,j}-\delta_j ) \right) , \end{split} \]
where $W_{i,j}=\sum_{k \in [1,M] \cup H_1 \cup J_1 \cup \cdots \cup H_{i-1} \cup J_{i-1} \cup H_i \cup J_{j-1}} X_k$. Let $\mathcal{G}_{i,j}$ denote the $\sigma$-algebra generated by the variables $X_k$, $k \in [1,M] \cup H_1 \cup J_1 \cup \cdots \cup H_{i-1} \cup J_{i-1} \cup H_i \cup J_{j-1}$ and $\xi_i$, $\xi_j$. Note that $W_{i,j}$, $\delta_i$, $\delta_j$ and $S_k$, $k \in H_i$ are $\mathcal{G}_{i,j}$-measurable. Moreover, $S_{\ell}-W_{i,j}$ is the sum of $\ell -M-|H_1|-|J_1|-\cdots - |H_{i-1}|-|J_{i-1}|-|H_i|-|J_{j-1}|$ independent copies of $X_1$, and it is independent of $\mathcal{G}_{i,j}$. We can thus apply \eqref{maintool} to $X=S_{\ell}-W_{i,j}$ and a random variable $Y$ uniformly distributed on $[0,1)$, independent of $\mathcal{G}_{i,j}$ to obtain
\begin{multline}\label{conditionalgij}
\left| \E \left( (T_i-T_i^*)(T_j-T_j^*) \mid \mathcal{G}_{i,j} \right) \right| \ll \\ \sum_{\substack{k \in H_i \\ \ell \in H_j}} \left| f(S_k)-f(S_k-\delta_i ) \right| \psi (\ell -M-|H_1|-|J_1|-\cdots - |H_{i-1}|-|J_{i-1}|-|H_i|-|J_{j-1}|) .
\end{multline}
Next, we condition with respect to the $\sigma$-algebra $\mathcal{F}_i$ generated by $X_k$, $k \in J_{i-1}$ and $\xi_i$. Let $Y_i=\sum_{k \in J_{i-1}}X_k$, and let us write $|f(S_k)-f(S_k-\delta_i )| = |f((S_k-Y_i)+Y_i)-$ $f((S_k-Y_i)+Y_i-\delta_i)|$, as before. Note that $Y_i$ and $\delta_i$ are $\mathcal{F}_i$-measurable, while $S_k-Y_i$ is independent of $\mathcal{F}_i$. Let us apply \eqref{maintool} to $X=S_k-Y_i$, a random variable $Y$ uniformly distributed on $[0,1)$, independent of $\mathcal{F}_i$ and $g(x)=|f(x+Y_i)-$ $f(x+Y_i-\delta_i)|$. As observed before, $\| g \|_1 \ll |\delta_i| \le \psi (|J_{i-1}|)$, hence
\begin{equation}\label{conditionalfi}
\E \left( |f(S_k)-f(S_k-\delta_i)| \mid \mathcal{F}_i \right) \ll \psi (|J_{i-1}|) + \psi (k-|J_{i-1}|) \ll i^{-c-\delta /2} .
\end{equation}
Taking the (total) expectation in \eqref{conditionalgij} and \eqref{conditionalfi}, and summing over $k \in H_i$ we finally deduce
\begin{multline}\label{covarianceestimate}
\left| \E (T_i-T_i^*)(T_j-T_j^*) \right| \ll \\ i^{-\delta /2} \sum_{\ell \in H_j} \psi (\ell -M- |H_1|-|J_1|-\cdots - |H_{i-1}|-|J_{i-1}|-|H_i|-|J_{j-1}|) .
\end{multline}
Now fix $2 \le i <R$. In the case $j=i+1$ the sum on the right hand side of \eqref{covarianceestimate} is $\ll 1$, thus $|\E (T_i-T_i^*)(T_{i+1}-T_{i+1}^*)| \ll i^{-\delta /2}$. In the case $j-i \ge 2$ we have $(\ell -M- |H_1|-|J_1|-\cdots - |H_{i-1}|-|J_{i-1}|-|H_i|-|J_{j-1}|) \gg$ $(j-i-1) j^c$, and hence $\left| \E (T_i-T_i^*)(T_j-T_j^*) \right| \ll (j-i-1)^{-1-\delta} j^{-c \delta} i^{-\delta /2}$. Therefore $\sum_{j=i+1}^{R} \left| \E (T_i-T_i^*)(T_j-T_j^*) \right| \ll i^{-\delta /2}$ for any fixed $2 \le i <R$. Summing over $2 \le i <R$ we finally obtain that the contribution of the off-diagonal terms in \eqref{expandsq} is also $\ll R^{1-\delta /2}$.
\end{proof}

\section{Exponential bounds}

The crucial step in the proof of laws of the iterated logarithm is to prove sharp bounds for the tail probabilities of the partial sums of the considered random variables. In this section we prove such an inequality for shifted sums $\sum_{k=M+1}^{M+N} f(S_k)$.
\begin{lem}\label{expbound}
Assume that $\psi (k) \ll k^{-(1+\delta)}$ for some $0<\delta<1$, and fix a constant $0<\rho <\delta /16$. For any integers $M \ge 0$ and $N \gg 1$, any $f \in \mathcal{F}$ such that $V(f) \le 2$ and $\| f \|_1 \gg N^{-1/2-\delta/8}$, and any $t \gg 1$ we have
\[ \P \left( \left| \sum_{k=M+1}^{M+N} f(S_k) \right|  \geq t \| f\|_1^\rho \sqrt{N} \right) \ll \exp \left( -C \| f \|_1^{-\rho} t_0 t \right) + t^{-2} \| f \|_1^{-2\rho} N^{-1/2-\delta /8} \]
with some constant $C>0$ and $t_0=\min \{ t, N^{\delta /16} \}$. The implied constants and $C$ depend only on the distribution of $X_1$, $\delta$ and $\rho$.
\end{lem}

\begin{proof} Let $c=1-\delta/2$ and $c'=(c+\delta/2)/(1+\delta)$. Let us split the set $\{ M+1, M+2, \dots \}$ into blocks $H_1, J_1, H_2, J_2, \ldots$ such that $|H_i|=\lceil i^c \rceil$ and $|J_i|=\lceil i^{c'} \rceil$ for all $i \ge 1$. Let $T_i, D_i$, $i \ge 1$ and $T_i^*, D_i^*$, $i \ge 2$ be as in Section \ref{approximationsubsection}. Fix a large enough integer $N$, and let $R$ be such that $M+N \in H_R \cup J_R$. Clearly $R = \Theta(N^{1/(c+1)})$, and $|T_i| \ll R^c \ll N^{(2-\delta)/(4-\delta)}$ for all $1 \le i \le R $, and the same holds for $D_i$. Since $N^{(2-\delta)/(4-\delta)} \ll \| f \|_1^{\rho} N^{1/2}$, the contribution of any individual block to $\sum_{k=M+1}^{M+N} f(S_k)$ is negligible. Therefore we may assume that $M+N=\max J_R$, and it will be enough to prove the claim for $\sum_{i=2}^R (T_i+D_i)$ instead of $\sum_{k=M+1}^{M+N} f(S_k)$. Recall that $T_i^*$, $i=2,3,\dots, R$ are independent, and $\E T_i^* =0$; similarly, $D_i^*$, $i=2,3,\dots, R$ are independent, and $\E D_i^* =0$. By Lemma \ref{approxerror1} we have
\[ \E \left(\sum_{i=2}^R (T_i-T_i^*)\right)^2 \ll R^{1-\delta/2} \ll N^{(1-\delta/2)/ (2-\delta /2)} \le N^{1/2-\delta /8}, \]
and the same holds for $\sum_{i=2}^R (D_i-D_i^*)$. Hence by the Chebyshev inequality
\begin{equation}\label{chebyshevineq}
\P \left( \left| \sum_{i=2}^R (T_i+D_i) -\sum_{i=2}^R (T_i^*+D_i^*)\right| \ge t \| f \|_1^\rho \sqrt{N} \right) \ll t^{-2} \| f \|_1^{-2\rho} N^{-1/2-\delta /8}.
\end{equation}
From \eqref{ti*di*variance} and $C(f,f) \ll \| f \|_1$ we get
\[ \sum_{i=2}^R \E T_i^{*2} \ll \| f \|_1 N + \sum_{i=2}^R i^{c/(1+\delta )} \ll \| f \|_1 N + N^{1-\delta /8} \ll \| f \|_1^{3 \rho} N . \]
Assume $|\lambda |<\left( 4 \lceil R^c \rceil \right)^{-1}$, then $|\lambda T^*_i|
\leq 2|\lambda| |H_i| \le 1/2$, and thus using the previous relations and $e^x \leq 1 + x + x^2$ for $|x| \leq 1/2$  we get
\begin{align*}
& \E \exp \left( \lambda \sum_{i=2}^R  T^*_i \right)
 = \prod_{i=2}^R \E e^{\lambda T^*_i}
 \leq \prod_{i=2}^R  \E  (1 + \lambda T^*_i + \lambda^2
T^{*2}_i) \nonumber \\
& = \prod_{i=2}^R (1 + \lambda^2 \E  {T^{*2}_i})
\leq \exp \left(\lambda^2 \sum_{i=2}^R \E T_i^{*2}\right) \leq \exp( C_1 \lambda^2 \| f \|_1^{3 \rho} N)
\end{align*}
with some constant $C_1>0$. Note that $\lambda = C_2 t_0 \| f \|_1^{-2\rho} N^{-1/2} \le C_2 \| f \|_1^{-2 \rho} N^{-1/2+\delta/16}$ with a small enough constant $C_2>0$ satisfies $|\lambda | < \left( 4 \lceil R^c \rceil \right)^{-1}$. Using the Markov inequality we thus get
\[ \begin{split} \P \left( \sum_{i=2}^R T^*_i \geq t \| f \|_1^\rho \sqrt{N} \right) &\le \exp \left(- \lambda t \| f \|_1^\rho \sqrt{N}\right) \E \exp \left(\lambda \sum_{i=2}^R T_i^* \right) \\
&\le \exp\left(- \lambda t \| f \|_1^\rho \sqrt{N} +C_1 \lambda^2 \| f \|_1^{3 \rho} N \right)   \\
&\le \exp\left( -C_2 \| f \|_1^{-\rho} t_0 t  +C_1 C_2^2 \| f \|_1^{-\rho} t_0^2 \right) \\
&\le \exp\left( -C  \| f \|_1^{-\rho} t_0 t \right) \end{split} \]
with some constant $C>0$. Using $-\lambda$ instead, we get that the same upper estimate holds for $\P \left( \sum_{i=2}^R T_i^* \le -t \| f \|_1^{\rho} \sqrt{N} \right)$. Repeating the same arguments for $\sum_{i=2}^R D_i^*$, we altogether have
\[ \P \left( \left| \sum_{i=2}^R (T_i^*+D_i^*) \right| \ge t \| f \|_1^{\rho} \sqrt{N} \right) \le 4 \exp (-C \| f \|_1^{-\rho} t_0 t) . \]
The last relation, together with (\ref{chebyshevineq}) implies Lemma \ref{expbound}.
\end{proof}

\section{Chaining}

Recall that for any interval $J \subseteq [0,1)$ of length $\lambda (J)$, $f_J (x)=I_J(\{ x \})-\lambda (J)$ denotes the centered indicator of $J$, extended with period $1$. Lemma \ref{expbound} in the previous section gives a useful estimate for the tail probabilities
$$\P \left( \left| \sum_{k=M+1}^{M+N} f_J(S_k) \right|  \geq   t \lambda(J)^\rho \sqrt{N} \right)$$
 which leads to the upper half of the LIL for the partial sums $\sum_{k=1}^N f_J(S_k)$ for a fixed interval $J \subseteq [0,1)$ and along an exponentially growing subsequence of $N$'s. Proving the LIL for the discrepancy $D_N(S_k)$ requires, on the other hand, tail estimates for the variable
$$ \max_{2^n\le N < 2^{n+1}} \sup_{J \subseteq [0,1)}\left|\sum_{k=1}^N f_J(S_k)\right|.$$
Such inequalities will be proved in this section, using the method of chaining. We follow closely Philipp \cite{PH}.

\begin{lem}\label{chaining1}
Assume that $\psi (k) \ll k^{-(1+\delta)}$ for some $0<\delta <1$. For any $a \in [0,1)$, any $0<\Delta \le 1-a$ and any integer $N \gg 1$ we have
\[ \P \left( \sup_{0<y\le \Delta} \left|\sum_{k=1}^N f_{[a,a+y)} (S_k) \right|\gg \Delta^{\delta /30}\sqrt{N}+ N^{1/2-\delta /30}\right)\ll \Delta^4+ N^{-\delta /30}. \]
The implied constants depend only on the distribution of $X_1$ and $\delta$.
\end{lem}

\begin{proof} For any $h \geq 1$, $1 \leq j \leq 2^h$ let $\varphi^{(j)}_h$ denote the indicator function of the interval $[(j - 1)2^{-h}, j2^{-h})$ extended with period $1$. We observe that if $0<y \le \Delta$ has the dyadic expansion
\begin{equation}\label{dyadic}
y = \sum^\infty_{\ell = 1} \e_\ell 2^{-\ell}, \qquad \e_\ell = 0,1,
\end{equation}
and $H \geq 1$ is an arbitrary integer, then
\begin{equation}\label{2}
\sum^H_{h=1} \varrho_h(x) \leq I_{[0,y)}(\{ x\} ) \leq \sum^H_{h
= 1} \varrho_h(x) + \sigma_H(x)
\end{equation}
where $\varrho_h$ is the indicator function of
$\Big[\sum\limits^{h-1}_{\ell = 1} \e_\ell 2^{-\ell},
\sum\limits^h_{\ell=1} \e_\ell 2^{-\ell} \Big)$ and
$\sigma_H$ is the indicator function of $\Big[
\sum\limits^H_{\ell=1} \e_\ell 2^{-\ell}, \sum\limits^H_{\ell=1}
\e_\ell 2^{-\ell} + 2^{-H}\Big)$, both extended with period $1$. For $\e_h = 0$
clearly $\varrho_h \equiv 0$, and thus (\ref{2}) remains valid if in the
sums we keep only those terms where $\e_h=1$. Also, for $\varepsilon_h = 1$,
$\varrho_h$ coincides with one of the $\varphi^{(j)}_h$, and $\sigma_H$ also
coincides with one of the $\varphi^{(j)}_H$. It follows that
\[ \begin{split} f_{[0,y)}(x) &\leq \underset{1\le h\le H}{\sum\nolimits^*}
(\varrho_h(x) - 2^{-h}) + (\sigma_H(x)-2^{-H}) +2^{-H}, \\ f_{[0,y)}(x) &\ge \underset{1\le h\le H}{\sum\nolimits^*}  (\varrho_h(x) - 2^{-h}) -2^{-H}, \end{split} \]
where $\sum^*$ means that the summation is extended only for those $h$ such that $\varepsilon_h=1$. Setting $x=S_k-a$ and summing for $1\le k\le N$ we get
\[ \begin{split} \sum_{k=1}^N f_{[a,a+y)}(S_k) &\leq \sum_{k=1}^N \underset{1\le h\le H}{\sum\nolimits^*}
(\varrho_h (S_k-a)- 2^{-h}) + \sum_{k=1}^N (\sigma_H (S_k-a)-2^{-H}) +N 2^{-H}, \\ \sum_{k=1}^N f_{[a,a+y)}(S_k) &\ge \sum_{k=1}^N \underset{1\le h\le H}{\sum\nolimits^*}  (\varrho_h (S_k-a)- 2^{-h}) -N 2^{-H}. \end{split} \]
Letting $A= \log_2 \Delta^{-1}$, observe that the dyadic expansion (\ref{dyadic}) of any $0<y\le \Delta$ starts with an index $\ell\ge A$. Thus it follows that for any integers $N \ge 1$ and $H \ge 1$
\begin{equation} \label{3}
\sup_{0<y \le \Delta} \left| \sum_{k=1}^N f_{[a,a+y)} (S_k) \right|
\le 2  \sum_{A \le h \le H} \max_{1 \le j \le 2^h} \left| \sum_{k=1}^N \left( \varphi^{(j)}_h
(S_k-a) - 2^{-h} \right) \right| + N2^{-H}.
\end{equation}
Choose $H=\lfloor (1/2+\delta /30) \log_2 N \rfloor$, and introduce the events
\begin{align*}
G(N,j,h) &= \left\{ \left| \sum_{k=1}^N \left( \varphi_h^{(j)} (S_k-a) -2^{-h} \right) \right| \geq D 2^{-(\delta /30) h} \sqrt{N}  \right\},\\
G_N &= \bigcup_{A\le h \leq H} \bigcup_{1 \le j \leq 2^h} G(N,j,h) , \nonumber
\end{align*}
where $D>0$ is a large constant, to be chosen. For any $h\le H$ we have $2^{-h} \ge N^{-1/2-\delta /30}$, and thus applying Lemma \ref{expbound} to the function \mbox{$\varphi^{(j)}_h (x-a)-2^{-h}$} with $t=t_0=D$ and $\rho=\delta /30$ we get
\[ \P (G(N, j, h)) \ll \exp (-C D^2 2^{(\delta /30) h} ) + 2^{(\delta /15) h} N^{-1/2-\delta /8}, \]
and consequently
\begin{equation}\label{4}
\P (G_N) \ll \sum_{h\ge A} \exp( h -C D^2 2^{(\delta /30) h}) +
N^{-1/2-\delta /8} \sum_{h \leq (1/2+\delta /30 ) \log_2 N} 2^{(1+\delta /15 ) h}.
\end{equation}
Here $h-CD^2 2^{(\delta /30)h} \le -4h$ provided $D$ is large enough. The first sum in \eqref{4} is thus $\ll \exp (-4A) \ll \Delta^4$. The second term on the right hand side of (\ref{4}) is  $\ll N^{-\delta /30}$, hence $\P(G_N) \ll \Delta^4 + N^{-\delta /30}$. It follows that with the exception of a set with probability $\ll \Delta^4 + N^{-\delta /30}$, the right hand side of (\ref{3}) is
\[ \ll \sum_{A\le h\le H}  2^{-(\delta /30) h} \sqrt{N}+N 2^{-H} \ll 2^{-(\delta /30)A} \sqrt{N}+ N^{1/2-\delta /30}\ll \Delta^{\delta /30} \sqrt{N}+ N^{1/2-\delta /30} . \]
\end{proof}

\begin{lem}\label{chaining2}
Assume that $\psi (k) \ll k^{-(1+\delta)}$ for some $0<\delta <1$. For any $a \in [0,1)$, any $0<\Delta \le 1-a$ and any integer $n \gg 1$ we have
\[ \P \left( \max_{2^n\le N<2^{n+1}} \sup_{0< y\le \Delta} \left|\sum_{k=1}^N f_{[a,a+y)} (S_k) \right| \gg \Delta^{\delta /30}\sqrt{2^n\log\log 2^n} + \sqrt{2^n} \right) \ll n^{-2} . \]
The implied constants depend only on the distribution of $X_1$ and $\delta$.
\end{lem}

\begin{proof}
As in the proof of Lemma \ref{chaining1}, for any $h \geq 1$, $1 \leq j \leq 2^h$ let $\varphi^{(j)}_h$ denote the
indicator function of the interval $[(j - 1)2^{-h}, j2^{-h})$ extended with period $1$. For any integers $M\ge 0$ and $N\ge 1$ let
\[ F(M, N, j, h) = \left| \sum^{M+N}_{k=M+1} \left( \varphi^{(j)}_h (S_k-a) - 2^{-h} \right) \right| . \]
Following the proof of Lemma \ref{chaining1}, we get that for any integers $M \ge 0$, $N \ge 1$ and $H \ge 1$
\begin{equation}\label{gaSk}
\sup_{0<y \le \Delta} \left| \sum_{k=M+1}^{M+N} f_{[a,a+y)}(S_k) \right| \le 2 \sum_{A \le h \le H} \max_{1 \le j \le 2^h} F(M, N, j, h) + N 2^{-H},
\end{equation}
where $A=\log_2 \Delta^{-1}$. Let now $n\ge 1$, and let $2^n \leq N < 2^{n+1}$ have dyadic expansion $N=2^n+\beta_{n-1}2^{n-1}+\cdots +\beta_1 2 +\beta_0$, $\beta_i=0,1$. Letting $M_{\ell}=2^n+\beta_{n-1}2^{n-1}+\cdots +\beta_{\ell} 2^{\ell}$ we can decompose the interval of integers $[1,N]$ as
\[ [1,N] = [1,2^n] \cup \bigcup_{\substack{1 \le \ell \le n \\ \beta_{\ell -1}=1}} [M_{\ell}+1, M_{\ell}+2^{\ell -1}] . \]
Note that the union of $[M_{\ell}+1,M_{\ell}+2^{\ell -1}]$ over $1 \le \ell <n/2$ has total size $\ll 2^{n/2}$. Further, $M_{\ell}=2^n+m_{\ell}2^{\ell}$ for some integer $0 \le m_{\ell}<2^{n-\ell}$. Applying \eqref{gaSk} on the interval $[1,2^n]$ with $H=n/2$ and on the intervals $[M_{\ell}+1,M_{\ell} +2^{\ell -1}]$ for all $n/2 \le \ell \le n$ such that $\beta_{\ell -1}=1$ with the choice $H=2^{\ell /2}$, we get that
\[ \begin{split} \sup_{0<y \le \Delta} \left| \sum_{k=1}^N f_{[a,a+y)}(S_k) \right| \ll &\sum_{A \le h \le n/2} \max_{1 \le j \le 2^h} F(0, 2^n, j, h) + 2^{n/2} \\ + &\sum_{\substack{n/2 \le \ell \le n \\ \beta_{\ell -1}=1}} \left( \sum_{A \le h \le \ell /2} \max_{1 \le j \le 2^h} F(2^n+m_{\ell}2^{\ell}, 2^{\ell -1}, j, h) + 2^{\ell /2} \right) . \end{split} \]
Here $\sum_{n/2 \le \ell \le n} 2^{\ell /2} \ll 2^{n/2}$, therefore we altogether have
\begin{equation}\label{mainformula}
\begin{split} \max_{2^n \le N<2^{n+1}} \sup_{0<y \le \Delta} &\left| \sum_{k=1}^{N} f_{[a,a+y)} (S_k) \right| \ll \\ &\sum_{A \le h \le n/2} \max_{1 \le j \le 2^h} F(0, 2^n, j, h) + 2^{n/2} \\ &+ \sum_{n/2 \le \ell \le n} \sum_{A\le h\le  \ell /2} \max_{0 \le m < 2^{n-\ell}} \max_{1 \le j \le 2^h} F(2^n + m 2^\ell, 2^{\ell - 1}, j, h). \end{split}
\end{equation}
Let us introduce the events
\[ \begin{split}
G(n,j,h) &=  \left\{ F(0,2^n, j, h) \geq D 2^{-(\delta /30)h} \sqrt{2^n \log\log 2^n} \right\}, \\  G_n &= \bigcup_{A \le h \leq n/2} \bigcup_{1 \le j \leq 2^h} G(n,j,h), \\
G^*(n, j, h, m, \ell ) &= \left\{ F(2^n + m2^\ell, 2^{\ell - 1}, j,
h) \geq  D 2^{-(\delta /30)(h+n-\ell )} \sqrt{2^n \log\log 2^n} \right\}, \\
G_n^* &= \bigcup_{n/2 \leq \ell \leq n} \bigcup_{A\le h \leq \ell /2} \bigcup_{0 \le m < 2^{n -
\ell}} \bigcup_{1 \le j \leq 2^h} G^*(n,j,h,m,\ell ) ,
\end{split} \]
where $D>0$ is a large constant, to be chosen. Applying Lemma \ref{expbound} to the function $\varphi_h^{(j)}(x-a)-2^{-h}$ with $M = 0$, $N = 2^n$, $t=t_0= D \sqrt{\log\log 2^n}$ and $\rho=\delta /30$ we get
\[ \P (G(n, j, h)) \ll \exp (-C D^2 2^{(\delta /30)h} \log \log 2^n ) + 2^{(\delta /15)h} 2^{(-1/2-\delta /8)n}, \]
and consequently
\begin{equation} \label{gn}
\P (G_n) \ll \sum_{h=0}^{\infty} \exp( h-C D^2 2^{(\delta /30)h} \log \log 2^n) + 2^{(-1/2-\delta /8)n} \sum_{h \leq n/2} 2^{(1+\delta /15)h}.
\end{equation}
Here $h-CD^2 2^{(\delta /30) h} \log \log 2^n \le -h-2 \log \log 2^n$ provided $D$ is large enough, therefore the first sum in \eqref{gn} is $\ll \exp (-2 \log \log 2^n ) \ll n^{-2}$. The second term in \eqref{gn} is $\ll 2^{-(\delta /30)n}$, hence $\P (G_n) \ll n^{-2}$.

Similarly, applying Lemma \ref{expbound} to the function $\varphi_h^{(j)}(x-a)-2^{-h}$ with $M = 2^n + m2^\ell$, $N = 2^{\ell
- 1}$, $t = D 2^{(1/2-\delta /30)(n-\ell )+1/2} \sqrt{\log \log 2^n}$, $t_0 \ge D \sqrt{ \log \log 2^n}$ and $\rho = \delta /30$ we get
\[ \begin{split} \P (G^* (n,j,h,m,\ell )) \ll &\exp (-C D^2 2^{(\delta /30)h} 2^{(1/2-\delta /30)(n-\ell )} \log \log 2^n) \\&+ 2^{(-1+\delta /15) (n-\ell )} 2^{(\delta /15) h} 2^{(-1/2-\delta /8)\ell}, \end{split} \]
and consequently
\begin{equation}\label{gn*}
\begin{split} \P (G_n^*)\ll &\sum_{n/2 \leq \ell \leq n} \sum_{h=0}^{\infty} \exp (h+n-\ell -C D^2 2^{(\delta /30)h} 2^{(1/2-\delta /30)(n-\ell )} \log\log 2^n) \\
&+ \sum_{n/2 \leq \ell \leq n} \sum_{h \leq \ell /2}  2^h 2^{n-\ell} 2^{(-1+\delta /15)(n-\ell )} 2^{(\delta /15) h} 2^{(-1/2-\delta /8) \ell} . \end{split}
\end{equation}
Here $h+n-\ell -C D^2 2^{(\delta /30)h} 2^{(1/2-\delta /30)(n-\ell )} \log \log 2^n \le -h-(n-\ell) -2 \log \log 2^n$ provided $D$ is large enough, therefore the first sum in \eqref{gn*} is $\ll \exp (-2\log \log 2^n) \ll n^{-2}$. Since the second sum in \eqref{gn*} is $\ll 2^{-(\delta /80)n} \ll n^{-2}$, we have $\P (G_n^*) \ll n^{-2}$.

It follows that with the exception of a set with probability $\ll n^{-2}$, the right hand side of \eqref{mainformula} is
\[ \begin{split} &\ll \sum_{h \ge A} 2^{-(\delta /30)h} \sqrt{2^n \log \log 2^n} +2^{n/2} + \sum_{n/2 \le \ell \le n} \sum_{h \ge A} 2^{-(\delta /30)(h+n-\ell )} \sqrt{2^n \log \log 2^n} \\ &\ll \Delta^{\delta /30} \sqrt{2^n \log \log 2^n} +2^{n/2} . \end{split} \]
\end{proof}

\section{Proof of Theorem \ref{theorem2}}

In Lemma \ref{varianceprop1} we proved that $C(f,g)$ is bilinear, symmetric and positive semidefinite in $f,g \in \mathcal{F}$. In particular, $\Gamma (s,t)$ is symmetric and positive semidefinite. From the triangle inequality we get that
\[ \begin{split} |\E (f_{[0,s)}(U) - f_{[0,s')}(U)) f_{[0,t)}(U+S_k)| &\le 2 |s-s'|, \\ |\E f_{[0,s)}(U)(f_{[0,t)}(U+S_k) - f_{[0,t')}(U+S_k))| &\le 2 |t-t'| \end{split} \]
for any $s,s',t,t' \in [0,1]$ and any $k \ge 1$. Hence $\E f_{[0,s)}(U) f_{[0,t)}(U+S_k)$ is continuous on the unit square for any fixed $k \ge 1$. As observed in the proof of Lemma \ref{varianceprop1} (i), each term in the series defining $\Gamma (s,t)$ has absolute value at most $\psi (k) \ll k^{-(1+\delta)}$, and so the continuity of $\Gamma (s,t)$ follows.
\begin{lem}\label{multclt}
For any fixed $0\le t_1<t_2<\cdots <t_r \le 1$ the random vector ${\bf Y}_k=(f_{[0, t_1)} (S_k), f_{[0, t_2)} (S_k), \ldots, f_{[0, t_r)} (S_k))$ satisfies
\[ N^{-1/2} \sum_{k=1}^N {\bf Y}_k \overset{d}{\longrightarrow} \mathcal{N}(\bf{0}, \Sigma) , \]
where $\mathcal{N}(\bf{0}, \Sigma)$ denotes the multidimensional normal distribution with mean $\bf{0}$ and covariance matrix $\mathbf{\Sigma}= (\Gamma (t_i, t_j))_{1\le i, j\le r}$.
\end{lem}

\begin{proof} Let $c=1-\delta /8$ and $c'=(c+\delta /2)/(1+\delta)$. Note that $c'<c$. Let us decompose the set of positive integers into consecutive blocks $H_1, J_1, H_2, J_2, \dots$ such that $|H_i|=\lceil i^c \rceil$ and $|J_i|=\lceil i^{c'} \rceil$ for all $i \ge 1$. Put $f_j=f_{[0,t_j)}$, $j=1,2,\dots, r$. Following the notation of Section \ref{approximationsubsection} let ${\bf T}_i=\left( T_i^{(f_1)}, T_i^{(f_2)}, \dots, T_i^{(f_r)} \right)$ and ${\bf T}_i^*=\left( T_i^{(f_1)*}, T_i^{(f_2)*}, \dots, T_i^{(f_r)*} \right)$, and let ${\bf D}_i$ and ${\bf D}_i^*$ be defined analogously. For any positive integer $N$ let $R=R(N)$ be such that $N \in H_R \cup J_R$. Note that $R=\Theta \left( N^{1/(c+1)} \right)$. Since $|{\bf T}_i| \ll |H_i| \ll R^c \ll N^{1/2-\delta /32}$ for all $1 \le i \le R$ and the same holds for $|{\bf D}_i|$, the contribution of any individual block is negligible, and we get
\begin{equation}\label{Ykfirstdecomposition}
\sum_{k=1}^N {\bf Y}_k = \sum_{i=2}^R ({\bf T}_i + {\bf D}_i) + O \left( N^{1/2-\delta /32} \right) .
\end{equation}
Applying Lemma \ref{approxerror1} to $f_1, f_2, \dots, f_r$ we get $\E \left| \sum_{i=2}^R \left( {\bf T}_i + {\bf D}_i - {\bf T}_i^* - {\bf D}_i^* \right) \right|^2 \ll R^{1-\delta /2}$, and hence from the Chebyshev inequality and the Borel--Cantelli lemma $\left| \sum_{i=2}^R \left( {\bf T}_i + {\bf D}_i - {\bf T}_i^* - {\bf D}_i^* \right) \right| \ll R^{1-\delta /8}$ a.s. Thus replacing ${\bf T}_i$ by ${\bf T}_i^*$ and ${\bf D}_i$ by ${\bf D}_i^*$, the right hand side of \eqref{Ykfirstdecomposition} changes by $\ll R^{1-\delta /8} \ll N^{1/2-\delta /32}$, and so
\begin{equation}\label{Ykseconddecomposition}
\sum_{k=1}^N {\bf Y}_k = \sum_{i=2}^R ({\bf T}_i^* + {\bf D}_i^*) + O \left( N^{1/2-\delta /32} \right) \quad \textrm{a.s.}
\end{equation}
Thus it will be enough to prove the multidimensional CLT for the right hand side of \eqref{Ykseconddecomposition}.

Recall that ${\bf T}_2^*, {\bf T}_3^*, \dots$ are independent, mean ${\bf 0}$ random vectors. From \eqref{ti*di*variance} we deduce
\[ \mathrm{Cov} \left( T_i^{(f_j)*} , T_i^{(f_k)*} \right) = \frac{1}{4} \left( \E T_i^{(f_j+f_k)*2} - \E T_i^{(f_j-f_k)*2} \right) \sim C(f_j,f_k) |H_i| . \]
Here $C(f_j,f_k)=\Gamma (t_j,t_k)$, thus letting ${\bf \Sigma}_i$ denote the covariance matrix of ${\bf T}_i^*$, we have ${\bf \Sigma}_i \sim |H_i| {\bf \Sigma}$. Since $\sum_{i=2}^R |H_i| \sim N$, we have $N^{-1}\sum_{i=2}^R {\bf \Sigma}_i \to {\bf \Sigma}$. Note that $|{\bf T}_i^*| \ll |H_i| \ll R^c \ll N^{1/2-\delta /32}$ for all $2 \le i \le R$. Thus ${\bf T}_2^*, {\bf T}_3^*, \dots$ satisfies the Lindeberg condition, and by the classical multidimensional CLT we have
\begin{equation}\label{ti*multidimCLT}
N^{-1/2} \sum_{i=2}^R {\bf T}_i^* \overset{d}{\longrightarrow} \mathcal{N}(\bf{0}, \Sigma).
\end{equation}
Similarly, ${\bf D}_2^*, {\bf D}_3^*, \dots$ are independent, mean ${\bf 0}$ random vectors. From \eqref{ti*di*variance} we deduce
\[ \sum_{i=2}^R \E |{\bf D}_i^*|^2 \ll \sum_{i=2}^R |J_i| \ll R^{1+c'} = o(N), \]
therefore $N^{-1/2} \sum_{i=2}^R {\bf D}_i^* \overset{\P}{\longrightarrow} 0$. The last relation together with \eqref{ti*multidimCLT} imply $N^{-1/2} \sum_{i=2}^R ({\bf T}_i^* + {\bf D}_i^*) \overset{d}{\longrightarrow} \mathcal{N}(\bf{0}, \Sigma)$, as claimed.
\end{proof}

Lemma \ref{multclt} implies the convergence of the finite dimensional distributions of the sequence $\sqrt{N} (F_N(t)-t)$ in (\ref{donsker}) to those of $K(t)$, and to prove Theorem \ref{theorem2} it remains to prove the tightness of the sequence in $D[0,1]$. To this end, fix $0<\e<1$ and let $m$ be the positive integer such that $1/(m+1) \le \e< 1/m$. Applying Lemma \ref{chaining1} with $a=j/m$, $j=0,1,\dots, m-2$ and $\Delta =2/m$, it follows that
\begin{multline*}
\P \left(\max_{0\le j \le m-2} \sup_{0< y \le 2/m} \left|\sum_{k=1}^N f_{[j/m,j/m+y)}(S_k) \right|\gg m^{-\delta /30}\sqrt{N}+  N^{1/2-\delta /30}\right) \\ \ll m^{-3}+m N^{-\delta /30}.
\end{multline*}
Note that for any $s,t \in [0,1]$ with $s<t \le s+\e$ we have $\sqrt{N}|(F_N(t)-t)-(F_N(s)-s)|$ $=N^{-1/2} \sum_{k=1}^N f_{[s,t)}(S_k)$, and that $[s,t)$ is a subset of $[j/m,(j+2)/m)$ for some $0 \le j \le m-2$. This implies that with the exception of a set with probability $\ll \e^3 + \e^{-1} N^{-\delta /30}$ the fluctuation of the process $\sqrt{N} (F_N(t)-t)$ over any subinterval of $[0,1]$ with length $\le \e$ is $ \ll \e^{\delta /30} + N^{-\delta /30}$. By Theorem 15.2 of Billingsley \cite[Chapter 3]{BIL}, the sequence $\sqrt{N} (F_N(t)-t)$ is tight in $D[0,1]$. This completes the proof of Theorem \ref{theorem2}.

\section{Proof of Theorem \ref{theorem3}}

\begin{lem}\label{multlil} Fix $0\le t_1<t_2<\cdots <t_r \le 1$, and consider the random vector ${\bf Y}_k=(f_{[0, t_1)} (S_k), f_{[0, t_2)} (S_k), \ldots, f_{[0, t_r)} (S_k))$. With probability $1$, the set of limit points of the sequence $(2N\log\log N)^{-1/2} \sum_{k=1}^N {\bf Y}_k$ is the ellipsoid
\begin{equation}\label{ellips}
\left\{ {\bf \Sigma x} \in \mathbb{R}^r : \langle {\bf x}, {\bf \Sigma x}  \rangle  \le 1 \right\}
\end{equation}
in the range $\left\{ {\bf \Sigma x} : {\bf x} \in \mathbb{R}^r \right\}$ of the matrix $\mathbf{\Sigma}= (\Gamma (t_i, t_j))_{1\le i, j\le r}$.
\end{lem}

\begin{proof} Lemma \ref{multlil} can be proved in the same way as Lemma \ref{multclt}. The only difference is that after deducing \eqref{Ykseconddecomposition}, instead of the multidimensional CLT we apply the multidimensional version of Kolmogorov's LIL due to Berning \cite{BER}. This way we get that with probability $1$, the set of limit points of the sequence $(2N \log \log N)^{-1/2} \sum_{i=2}^R {\bf T}_i^*$ is the ellipsoid \eqref{ellips}, and that $(2N \log \log N)^{-1/2} \sum_{i=2}^R {\bf D}_i^* \to {\bf 0}$ a.s.
\end{proof}

Using Lemma \ref{multlil} and Lemma \ref{chaining2}, the proof of Theorem \ref{theorem3} can be easily completed. Let
\[ \alpha_N (t)= \frac{N (F_N(t)-t)}{\sqrt{2 N \log \log N}}, \quad 0\le t \le 1, \]
and let $L$ denote the (random) class of limit functions of subsequences of $\alpha_N$ in the $D[0,1]$ topology. Similarly to the  proof of the functional CLT, which requires tightness and the convergence of the finite dimensional distributions, to prove the functional LIL in Theorem \ref{theorem3} we will show that with probability $1$ the sequence $\alpha_N$ is relatively compact in the $D[0,1]$ topology, and that for any fixed finite subset $T \subset [0,1]$ we have $\{ f|_T : f \in L \} = \{ f|_T : f \in B (\Gamma ) \}$ a.s.

To prove relative compactness, fix $0<\e <1$ and let $m$ be the positive integer such that $1/(m+1) \le \e <1/m$. Applying Lemma \ref{chaining2} with $a=j/m$, $j=0,1,\dots, m-2$ and $\Delta = 2/m$, it follows that
\begin{multline*}
\P \Bigg( \max_{0 \le j \le m-2} \max_{2^n \le N < 2^{n+1}} \sup_{0<y \le 2/m} \left|\sum_{k=1}^N f_{[j/m,j/m+y)}(S_k) \right| \\ \gg m^{-\delta /30} \sqrt{2^n \log \log 2^n} +\sqrt{2^n} \Bigg) \ll m/n^2 .
\end{multline*}
By the Borel--Cantelli lemma this implies that with probability $1$, for any $0<\e <1$ there exists a (random) integer $N_0$ such that for any $N\ge N_0$ the fluctuation of the process $\alpha_N (t)$ over any subinterval of $[0,1]$ of length $\le \e$ is $\ll \e^{\delta /30}$. Applying Lemma \ref{chaining2} with $a=0$ and $\Delta =1$, from the Borel--Cantelli lemma we similarly deduce that with probability $1$, $\sup_{t \in [0,1]} |\alpha_N (t)| \ll 1$. By Theorem 14.3 of Billingsley \cite[Chapter 3]{BIL}, with probability $1$ the sequence $\alpha_N(t)$ is relatively compact in $D[0,1]$. Further, it is easy to see from the fluctuation estimate that with probability $1$, any subsequence of $\alpha_N(t)$ convergent in $D[0,1]$ has a continuous limit function. Hence with probability $1$, $L$ is a compact subset of $C[0,1]$ and any subsequence of $\alpha_N$ convergent in $D[0,1]$ is also uniformly convergent.

Next, fix a finite set $T \subset [0,1]$. Let $H(\Gamma )$ and $H\left( \Gamma|_{T \times T} \right)$ denote the reproducing kernel Hilbert spaces determined by the covariance functions $\Gamma$ and $\Gamma|_{T \times T}$, respectively. Then $H \left( \Gamma|_{T \times T} \right) = \{ f|_T : f \in H (\Gamma) \}$; moreover, $\| g \|_{H \left( \Gamma|_{T \times T} \right)} = \min \{ \| f \|_{H(\Gamma)} : f|_T =g \}$, see Aronszajn \cite[Section I.5]{A}. In particular, the closed unit ball of $H\left( \Gamma|_{T \times T} \right)$ is $\{ f|_T : f \in B(\Gamma) \}$. Lemma \ref{multlil} thus shows that with probability $1$, the set of limit points of $\alpha_N|_T$ is identical with $\{ f|_T : f \in B( \Gamma) \}$. The relative compactness of $\alpha_N$ implies that in fact $\{ f|_T : f \in L \} = \{ f|_T : f \in B (\Gamma ) \}$ a.s., as claimed.

Now fix finite sets $T_1 \subset T_2 \subset \cdots \subset [0,1]$ such that $\bigcup_{n=1}^{\infty} T_n$ is dense in $[0,1]$. Observe that if $K$ and $K'$ are compact subsets of $C[0,1]$ such that $\{ f|_{T_n} : f \in K \} = \{ f|_{T_n} : f \in K' \}$ for all $n \ge 1$, then $K=K'$. By Lemma 3 of Oodaira \cite{OO}, the continuity of $\Gamma (s,t)$ implies that $B(\Gamma )$ is compact in $C[0,1]$. The arguments above show that with probability $1$, $L$ is a compact subset of $C[0,1]$ and $\{ f|_{T_n} : f \in L \} = \{ f|_{T_n} : f \in B(\Gamma ) \}$ for all $n \ge 1$. Therefore $L=B(\Gamma )$ a.s. This finishes the proof of Theorem \ref{theorem3}.

\section{Proof of Theorem \ref{theorem1}} The proof of Theorem \ref{theorem1} is similar to, but is actually much simpler than that of Theorems \ref{theorem2} and \ref{theorem3}. The claims on $C(f,g)$ were proved in Lemma \ref{varianceprop1}. If in the proof of Lemmas \ref{multclt} and \ref{multlil} we replace the vector-valued random variables ${\bf Y}_k$, ${\bf T}_i$, ${\bf T}_i^*$, ${\bf D}_i$ and ${\bf D}_i^*$ by the real-valued random variables $f(S_k)$, $T_i^{(f)}$, $T_i^{(f)*}$, $D_i^{(f)}$ and $D_i^{(f)*}$, respectively, where $f \in \mathcal{F}$ is fixed, we get that the sum $\sum_{k=1}^N f(S_k)$ satisfies the CLT and the LIL.

\end{document}